\documentclass[twoside,leqno]{article}

\usepackage[letterpaper]{geometry}

\usepackage{siamproceedings}

\usepackage[T1]{fontenc}
\usepackage{amsfonts}
\usepackage{graphicx}
\usepackage{epstopdf}
\usepackage{enumitem}
\usepackage{algorithmic}
\ifpdf
  \DeclareGraphicsExtensions{.eps,.pdf,.png,.jpg}
\else
  \DeclareGraphicsExtensions{.eps}
\fi

\newsiamremark{remark}{Remark}
\newsiamremark{hypothesis}{Hypothesis}
\crefname{hypothesis}{Hypothesis}{Hypotheses}
\newsiamthm{claim}{Claim}

\usepackage{amsmath,bbold,mathrsfs,amscd,lineno,hyperref}
\usepackage{mathtools,afterpage}
\usepackage{tikz}
\usepackage{mathrsfs}
\usepackage{yfonts,amsbsy} 

\usepackage{amsfonts,subcaption}
\usepackage{amssymb}
\usepackage{enumitem}

\usepackage{chngcntr}
\counterwithout{figure}{section}

\theoremstyle{plain}
\theoremheaderfont{\bfseries}
\theorembodyfont{\normalfont}
\renewtheorem{theorem}{Theorem}

\newtheorem{thm}{Theorem}[section]

\newtheorem{lem}[thm]{Lemma}

\newtheorem{prob}[thm]{Problem}
\newtheorem{Interprob}[thm]{Interpolation Problem}
\newtheorem{conj}[thm]{Conjecture}
\newtheorem{MLCconj}[thm]{MLC Conjecture}
\newtheorem{quest}[thm]{Question}

\DeclareFontFamily{U}{mathb}{\hyphenchar\font45}
\DeclareFontShape{U}{mathb}{m}{n}{ <5> <6> <7> <8> <9> <10> gen * mathb <10.95> mathb10 <12> <14.4> <17.28> <20.74> <24.88> mathb12 }{}
\DeclareSymbolFont{mathb}{U}{mathb}{m}{n}
\DeclareFontSubstitution{U}{mathb}{m}{n}
\DeclareMathSymbol{\selfmap}{3}{mathb}{"FD}

\newcommand{\C}{{\mathbb{C}}}
\newcommand{\wC}{{\widehat{\mathbb{C}}}}

\newcommand{\Q}{\mathbb{Q}}
\newcommand{\R}{\mathbb{R}}
\newcommand{\N}{\mathbb{N}}
\newcommand{\Z}{\mathbb{Z}}

\newcommand{\Disk}{\mathbb{D}}

\newcommand{\Width}{{\mathcal W}}

\usepackage{amsopn}

\newcommand{\codim}{{\operatorname{codim}}}
\newcommand{\wZ}{{\widehat Z}}

\newcommand{\diam}{\operatorname{diam}}

\renewcommand{\mod}{\operatorname{mod}}

\newcommand{\orb}{\operatorname{orb}}

\renewcommand{\sec}{{\operatorname{sec}}}

\newcommand{\hor}{{\operatorname{hor}}}
\newcommand{\ver}{{\operatorname{ver}}}
\newcommand{\tot}{{\operatorname{tot}}}
\newcommand{\hide}[1]{}

\newcommand{\Hors}{{\mathfrak H}}
\newcommand{\cyl}{{\operatorname{cyl}}}

\newcommand\MM{{\mathcal{M}}}
\newcommand\HH{{\mathcal{H}}}
\newcommand\Mol{\mathfrak M}

\newcommand\Jul{J}
\newcommand\filled{K}
\newcommand{\QL}{{\mathcal Q\mathcal L}}
\renewcommand{\bmod}{ {\operatorname{mod}^+}}

\newcommand{\bnd}{{\mathrm{bnd}}}

\newcommand{\ovlTheta}{{\overline \Theta}}

\makeatletter
\newcommand{\setword}[2]{%
  \phantomsection
  #1\def\@currentlabel{\unexpanded{#1}}\label{#2}%
}

\newcommand{\pp}{\mathfrak{p}}

\newcommand{\qq}{\mathfrak{q}}

\newcommand{\nn}{\mathfrak{n}}

\newcommand{\RR}{{\mathcal R}}

\title{On the MLC Conjecture and the \\Renormalization Theory in Complex Dynamics}

\author{Dzmitry Dudko\thanks{Department of Mathematics, Stony Brook University, USA}}

\begin{document}

\date{}

\maketitle

\begin{abstract}  In this Note, we present recent developments in the Renormalization Theory of quadratic polynomials and discuss their applications, with an emphasis on the MLC conjecture, the problem of local connectivity of the Mandelbrot set, and on its geometric counterparts. 
\end{abstract}

\section{Introduction.}\label{s:intro} 
\begin{figure}[t]  \centering
   \includegraphics[width=14cm,trim={0 15 0 15}, clip]{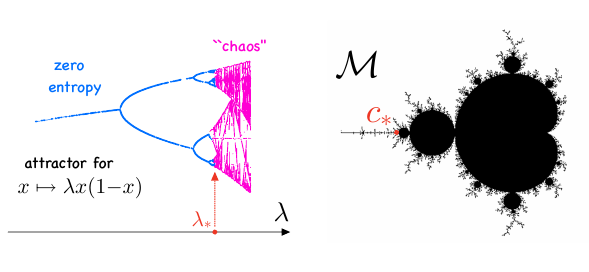}
   \caption{The Mandelbrot set (right) and the period-doubling cascade (left, colored in blue). The period-doubling Feigenbaum parameter $\lambda_* \simeq c_*$ marks the transition from zero entropy to positive entropy; it is where the ``chaos'' emerges. This transition is universal; its complexification is shown in Figure~\ref{Fig:Feigenb:C}.}
   \label{Fig:Mand_set}
\end{figure}

The Mandelbrot set $\MM$ (Figure~\ref{Fig:Mand_set}) is one of the most recognizable fractals in mathematics and beyond. It encodes the dynamical dependence of quadratic polynomials on the parameter. A central question in Complex Dynamics is the MLC Conjecture, also known as the \emph{Rigidity Conjecture}:
\begin{MLCconj}
    The Mandelbrot set $\MM$ is locally connected.
\end{MLCconj}

 This conjecture has various geometric and probabilistic counterparts; all of them are subjects of the Renormalization Theory of $\MM$, which analyzes \emph{first return maps} to small neighborhoods of special points. Nature and the classification of these first return maps is in focus of this Note. Let us start by outlining its context.

 In Section~\ref{s:intro} below, we give an informal introduction to renormalization ideas. Although the MLC is a topological statement, any substantial progress typically requires geometric input. Conjecture~\ref{cnj:hyp_comp} states a simple geometric property of hyperbolic components that easily implies the MLC. A foundation of any renormalization theory is \emph{a priori} bounds, see Conjecture \ref{cnj:unif_bounds}: as we zoom into the Julia set and endow it with the first return map, we expect precompactness. A similar precompactness is expected in the parameter plane, see~\S\ref{ss:resc:param}. In~\S\ref{ss:resc:dynam}, we classify (sequences of) first return maps as either \emph{renorm-expanding} or \emph{renorm-balanced}. The exposition in Section~\ref{s:intro} is conjectural, but it has been justified in many cases.

Historical background is presented in Section~\ref{s:History}, where we also sketch the evolution of key ideas and briefly comment on the higher-degree cases. The development of the renormalization theory can be broadly divided into several stages. First, as already mentioned, initial \emph{a priori} bounds are required. These are used in justifying various rigidity properties, e.g., the MLC. Next, a detailed geometric analysis of the convergence of first-return maps follows; a desired goal is to establish the \emph{hyperbolicity} of the associated operator. Finally, deep applications can be derived; for instance, the ``regular vs stochastic'' dichotomy describing the \emph{nature of chaos}. 


Section~\ref{s:RrnormTheory:MM} is dedicated to the theory of renormalization operators acting on infinite-dimensional spaces of maps. We recognize three types of renormalization regimes: \emph{near-Neutral}, if some hyperbolic components become unbounded in the respective (small) parameter neighborhoods relative to the scale, \emph{Quadratic-Like}, in the presence of intermediate hyperbolic components relative to the scale, and \emph{Puzzle} (it is intertwined with the renorm-expanding notion), if all hyperbolic components become negligible relative to the scale. We outline the quadratic-like (``ql'') renormalization theory in its pure, bounded-type form. Its hyperbolicity program is now complete by Theorem~\ref{thm:QL:hyperb}, thanks to Theorem~\ref{thm:MLC:Feigenb}. We summarize the status of the near-Neutral Case and mention how Theorem~\ref{thm:wZ} provides an opening towards the full hyperbolicity of the Neutral Renormalization and beyond, see~\S\ref{ss:MolecRenorm}. We will also discuss the puzzle regime in~\S\ref{ss:puzzleRegime} and put forward Problem~\ref{conj:hyperb:main} on whether all these three regimes can be neatly integrated.

In Section~\ref{s:NDR}, on the Near-Degenerate Regime, we illustrate the sources of \emph{a priori} bounds on dynamical scales representing definite positive-entropy ``$\ge \varepsilon>0$'' (primitive ql and puzzle cases that are ``$\varepsilon$-away'' from the main Molecule) and zero-entropy (Molecule: neutral and bounded-type satellite ql) settings. We conclude by summarizing the current status of the MLC as two problems: completing the remaining unbounded satellite cases and developing interpolation techniques between definite positive and zero entropy settings. The former relies on further development of the theory behind \emph{almost-invariant pseudo-Siegel disks}, Figure~\ref{Fig:NearDegRegimIllustr}. For the latter, we put forward a strategy based on partially invariant \emph{virtual Julia sets} which control only \emph{partial postcritical sets}. If these programs are achieved, they would establish \emph{a priori} bounds for \emph{all} quadratic polynomials and open a clear path towards the MLC (rigidity) and its geometric and probabilistic counterparts discussed in this Note.

\subsection{Renormalization: summary and motivation.} The renormalization $\RR_{\eta} p_c$ of a quadratic polynomial $p_c\colon z\mapsto z^2+c$ can be defined as the rescaled first return map to a small $\eta$-neighborhood of the critical value, see~\eqref{eq:dfn:Rf}.  A single iterate of $\RR_\eta p_c$ controls a long finite orbit of $p_c$. Hence, instead of iterating $p_c$, we wish to iterate renormalization
\begin{equation}
    \label{eq:renorm:f_c}
p_c\ \leadsto \ \RR_{\eta_1} p_c\ \leadsto\  \RR_{\eta_2} p_c\ \leadsto \ \dots\ \leadsto\ \RR_{\eta_n} p_c\ \leadsto\  \dots\qquad \overset{\text{``}\lim\text{''}}{\longrightarrow}\end{equation} 
and then deduce strong dynamical consequences for $p_c$. To realize such a vision, initial \emph{a priori} bounds (precompactness conditions) are required, see Conjecture~\ref{cnj:unif_bounds}. 

We have two cases. If the derivative $\big[\RR_{\eta_n} p_c\big]'(z)$ typically grows, then we should expect an almost linear phase-parameter relationship between $\MM$ and $\Jul_c$ at $c$ as illustrated in Figure~\ref{Fig:Exp_Case}. We call this regime \emph{renorm-expanding}, see~\S\ref{ss:resc:dynam}; it is controlled by \emph{puzzles}. On the other hand, if the derivative $\big[\RR_{\eta_n} p_c\big]'(z)$ is bounded, then there should be geometric limits for the renormalization, see~\eqref{eq:resc:limits}.  We refer to this regime as \emph{renorm-balanced}.

\subsection{Hyperbolic components.} 
\begin{figure}[t]  \centering
   \includegraphics[width=14cm,trim={0 22 0 20}, clip]{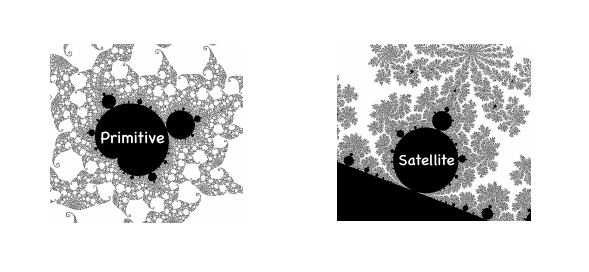}
   \caption{Primitive and satellite hyperbolic components. Around every hyperbolic component, there is an associated primitive or satellite little copy of the Mandelbrot set. It is expected, cf. Conjecture~\ref{cnj:hyp_comp}, that all hyperbolic components have uniform geometry. The difference between the satellite (zero entropy) and primitive (positive entropy) cases is profound and has many implications; e.g., see Figures~\ref{Fig:NDG:PositiveEntropy} and~\ref{Fig:NDR:ZeroEntropy}.
}
   \label{Fig:HypComp}
\end{figure}

For a hyperbolic component $\HH\subset \MM$, we denote by $\MM_\HH$ the associated small copy of $\MM$ centered at $\HH$. We recall that $\HH$ consists of parameters with an attracting periodic cycle and can be of two types, see Figure~\ref{Fig:HypComp}: satellite, when $\overline \HH$ is an analytic disk, and primitive, when $\overline \HH$ is an analytic image of the filled cardioid. On pictures, hyperbolic components always have a flawless shape.

\begin{conj}[Bounded geometry of $\HH\subset \MM_\HH$] \label{cnj:hyp_comp} Every closed hyperbolic component $\overline \HH$ is a uniform qc disk or a uniform qc filled cardioid; its diameter $\diam (\HH)$ is uniformly comparable to $\diam(\MM_\HH)$. 
\end{conj}
Conjecture~\ref{cnj:hyp_comp} easily implies the MLC: if $\MM_{\HH_n}$ is a sequence of strictly nested small copies of $\MM$, then the diameter of $\MM_{\HH_n}$ shrinks because of finiteness of the area:
\[\sum_n\big[\diam (\MM_{\HH_n})\big]^2\sim \sum_n\big[\diam (\HH_n)\big]^2\sim \sum_n\text{area}(\HH_n)<\infty,\qquad \qquad\text{thus }\bigcap_{n\ge 1} \MM_{\HH_n}\text{ is a singleton.}\]
This implies the MLC by Yoccoz's results, see~\eqref{eq:Yoccoz:results}. For primitive $\HH$, Conjecture~\ref{cnj:hyp_comp} is expected to follow from \emph{a priori} bounds of ql type \S\ref{ss:RR:QL:compact}; cf. Conjecture~\ref{cnj:unif_bounds} below. For satellite $\HH$, the hyperbolicity of the Molecule renormalization is likely necessary, see Conjecture~\ref{conj:hyp:Molecule}.

\subsection{Dynamical limits.}\label{ss:resc:dynam} Let us consider the dynamical plane of $p_c$ for a \emph{non-parabolic} $c\in \partial \MM$. (Renormalizations for parabolic $c\in \partial\MM$ arise through limits $c_n\to c$, see~\S\ref{ss:renorm:parablic}.) For a small $\eta>0$ and the associated $\eta$-neighborhood of the critical value $c$, the \emph{renormalization} $\RR_\eta p_c$ of $p_c$ is the \emph{rescaled first return map} of $p_c$ from this $\eta$-neighborhood back to itself; we write:
\begin{equation}
    \label{eq:dfn:Rf} \big[\RR_\eta p_c\big](z) \ =\ \Big[\big(p_c^{\qq}\big)^\text{$\eta$-rescaled}\ \Big] (z) \ =\ \frac{1}{\eta}\Big(p_c^{\qq} \big(\eta z+c\big)-c\Big),\qquad\qquad p_c(z)=z^2+c,
\end{equation}
where $\qq\equiv \qq(\eta z+c)$ is the first return time of the point $\eta z+c$ back to the $\eta$-neighborhood of $c$. The $\eta$-rescaling identifies the $\eta$-neighborhood of $c$ with the $1$-neighborhood of $0$.

\begin{figure}[t]  \centering
   \includegraphics[width=14cm,trim={0 10 0 8}, clip]{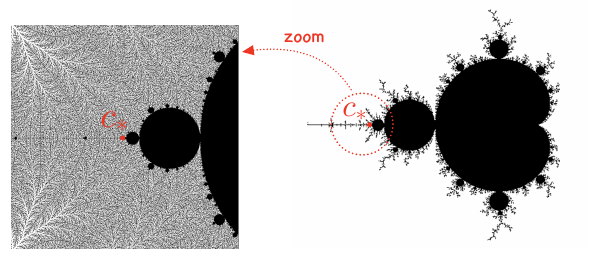}
   \caption{The universality of the Mandelbrot set at the period-doubling Feigenbaum parameter $c_*$, see also Figure~\ref{Fig:Mand_set}. The renormalization theory justifying this universality has been successfully developed in the 1990s, see~\S\ref{s:History}. The MLC at $c_*$ has been recently established in Theorem~\ref{thm:MLC:Feigenb}.}
   \label{Fig:Feigenb:C}
\end{figure}

\begin{conj}[Uniform \emph{a priori} bounds] \label{cnj:unif_bounds} With appropriate understanding, the family of renormalizations~\eqref{eq:dfn:Rf} is precompact uniformly over all $\eta>0$ and $c\in \partial \MM$. 
\end{conj}
\noindent In particular, as we zoom into the Julia set $\Jul_c$ of $p_c$ at $c$, the shape of  $\Jul_c$ (e.g., its veins) should stay within a precompact class of configurations. Let us detail more what we expect from Conjecture~\ref{cnj:unif_bounds}.

In~\eqref{eq:dfn:Rf}, different points can have different return times; and some points may not return to the selected neighborhood. Consider an iterate $p_c^\qq$ realizing a first return of some point in the $\eta$-neighborhood of $c$ back to the neighorhood. Rescaling  $p_c^\qq$, we obtain a \emph{globalized branch} $\RR_\eta^{(\qq)} p_c\equiv \big(p_c^{\qq}\big)^\text{$\eta$-rescaled} \colon \C\to \C$ of  $\RR_\eta p_c$. 

Consider now a sequence $\eta_n\to 0$ and the associated renormalizations $\RR_{\eta_n} p_c$. We refer to the $\eta_n$’s as \emph{dynamical scales}. Up to passing to a subsequence, we anticipate only two possibilities:

\begin{itemize}
\item[\setword{(Exp)}{item:Exp}] \emph{renorm-expanding} case: for all choices of branches, the derivative $\big[\RR^{(\qq_n)}_{\eta_n} p_c\big]'(z)$ tends to infinity for a typical (roughly, non-critical) $z\in \C$; or
\item [\setword{(Ren)}{item:Ren}] \emph{renorm-balanced} case: for some choices of globalized branches, there is a well-defined geometric limit, transcendental or polynomial (locally uniformly on $U_G$): \begin{equation}
    \label{eq:resc:limits} \exists\ \big[ G\colon U_G\to \C \big]\ =\ \lim_{n\to \infty} \big[\RR^{(\qq_n)}_{\eta_n} p_c\colon \C\to\C\big], \qquad \qquad  \eta_n\to 0.
\end{equation}
\end{itemize}

 \noindent We anticipate that the renorm-balanced case~\ref{item:Ren} occurs when certain hyperbolic components near $c$ do not shrink in the respective parameter scales (see~\S\ref{ss:resc:param}), so that the periodic cycles associated with such \emph{non-negligible} hyperbolic components prevent expansion; cf. Conjecture~\ref{cnj:hyp_comp}.

 By considering all sequences of dynamical scales $\eta_n\to 0$, we say that $c\in \MM$ is 
\begin{itemize}
    \item \emph{pure renorm-expanding,} if it is renorm-expanding with respect to any $\eta_n\to 0$; and
    \item \emph{pure renorm-balanced,} if it is renorm-balanced with respect to any $\eta_n\to 0$.
\end{itemize}


There are four main renormalization theories in the pure renorm-balanced case: primitive/satellite quadratic-like, near-parabolic, and near-Siegel; the last two are subtypes of near-neutral. These theories have many similarities as well as certain fundamental differences. The general case is a combination of all the above regimes, where various dynamical scales of $p_c$ can exhibit different renormalization flavors; see Section~\ref{s:RrnormTheory:MM} and Question~\ref{conj:hyperb:main}.

\subsection{Parameter limits.}\label{ss:resc:param} Recall~\cite{DH1} that the parameter and dynamical B\"ottcher functions provide conformal identification of $\wC\setminus \MM$ and $\wC\setminus \Jul_c$. Therefore, dynamical scales $\eta_n$ as in~\S\ref{ss:resc:dynam} have associated parameter scales $\delta_n$ defined up to appropriate choices. Let us write the parameter counterpart to~\eqref{eq:dfn:Rf}:
\begin{equation}
    \label{eq:dfn:R:MM} \RR_{(\delta,\eta)} \ \big\{p_{c+\tau}\big\}_{\tau\in \C}\ =\ \left(\big\{p_{c+\tau}\big\}_{\tau\in \C}\ \right)^\text{$(\delta, \eta)$-rescaled}\  \ =\ \frac{1}{\eta}\Big(p_{c+\delta \tau}^{\qq} \big(\eta z+c+\delta \tau\big)-c-\delta \tau\Big),
\end{equation}
i.e., all nearby maps $p_{c+\delta \tau}$ are $\eta$-rescaled with respect to the same return times as $p_c$, while the parameter dependence is $\delta$-rescaled. The notion of branches is defined accordingly.

Similar to Conjecture~\ref{cnj:unif_bounds}, we expect the precompactness of~\eqref{eq:dfn:R:MM} to hold uniformly over all $c$ and across all scales. In fact, we anticipate much stronger properties: convergence toward the limits should be exponentially fast, under an appropriate understanding, both in the parameter and dynamical planes, see~\S\ref{s:RrnormTheory:MM}. Such statements are commonly referred to as \emph{universality} and are illustrated in various figures throughout this Note.

\section{Historical Perspective.}\label{s:History} The problem of local connectivity of the Mandelbrot set $\MM$, known as the MLC Conjecture, was put forward in the 1980s by Douady and Hubbard, who also developed the foundational theory of $\MM$,~\cite{DH1, PolyLike}. It was the time when the first computer pictures of $\MM$ were generated; these pictures, ever since, have been providing invaluable guidance in research. 

Yoccoz' results from the early 1990s reduced the MLC Conjecture to the following rigidity statement~\cite{HY}:
\begin{multline}\label{eq:Yoccoz:results}
    \qquad \text{ if } \quad \MM_{1}\supsetneq\MM_2 \supsetneq \dots \supsetneq\MM_n \supsetneq\dots  \qquad\text{ are nested small copies of $\MM$,} \\ \text{then } \quad \bigcap_{n=1}^{\infty} \MM_n=\{c\}\quad\text{ is a singleton. }\qquad\end{multline}
This linked the MLC to the renormalization theory.

 Before the Mandelbrot set and its pictures, there has been the \emph{Logistic Family} which encodes the \emph{real vein} $\MM\cap \R$ of the Mandelbrot set:
\begin{equation}
\label{eq:log famil} x\mapsto \lambda x(1-x)\ \colon\  [0,1]\ \selfmap\ , \qquad\qquad \lambda\in [0,4].
\end{equation}
The discovery of \emph{chaos} is due to Poincar\'e (sensitive dependence combined with recurrence in the $3$-body problem, the late $19$th century). In an attempt to understand this phenomenon, much of the research activity in the second half of the 20th century was focused on the study of the Logistic Family (e.g., Ulam, Neumann in the 1940s, Lorenz, Myrberg, Sharkovsky in the 1960s, Misiurewicz, Guckenheimer in the 1970s). The Combinatorial Theory of~\eqref{eq:log famil} was designed by Milnor and Thurston in the late 1970s in terms of kneading sequences.

\begin{figure}[t]  \centering
   \includegraphics[width=14cm, trim={0 10 0 8}, clip]{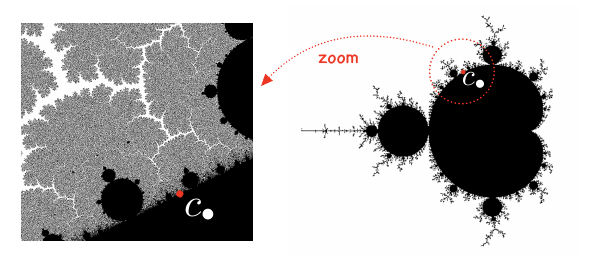}
   \caption{Universality of the Mandelbrot set at the golden-mean Siegel parameter $c_\bullet$. The dynamical universality of $z^2+c_\bullet$ is depicted in Figure~\ref{Fig:SiegelMaps}. Theorem~\ref{thm:hyperb:pacman} completed the program towards the hyperbolicity of the near-Siegel renormalization of bounded-type.}
   \label{Fig:SiegelFamil}
\end{figure}


Fifty years ago, the \emph{period-doubling universality} in~\eqref{eq:log famil} was discovered by Feigenbaum~\cite{F1} (the parameter universality) and, independently, by Tresser and Collet~\cite{CT} (the dynamical universality); see Figures~\ref{Fig:Mand_set} and~\ref{Fig:Feigenb:C}. This was the beginning of the Renormalization Theory in Dynamics; many renormalization ideas have been imported from Quantum Field Theory and Statistical Mechanics. The first computer-assisted proof of the period-doubling universality was presented by Lanford~\cite{Lan} in the early 1980s. A computer-free proof, with numerous extensions to other cases, was only achieved in the late 1990s using \emph{complex} methods.

In his address to the ICM in 1986 \cite{Sul:ICM86}, Sullivan proposed a program to approach the period-doubling universality based upon Teichm\"uller theory by employing the Douady-Hubbard quadratic-like theory. Sullivan partially realized this program in~\cite{Sullivan} by constructing a renormalization fixed point and justifying (non-exponential) convergence towards the fixed point within its hybrid class. The exponential convergence (the dynamical universality) was justified by McMullen in~\cite{McM3}. The program was completed by Lyubich in~\cite{LyuFeigen} who established the parameter universality.


Meanwhile, a \emph{complete} renormalization theory has been successfully developed for the Logistic Family~\cite{McM2,  LyuActa,GS}, culminating in the ``regular vs stochastic'' dichotomy \cite{LyuRegSt} by Lyubich; i.e., the Logistic Family became the first natural class of maps outside of the hyperbolic setting for which \emph{chaos was well understood}. Later, the dichotomy was refined to ``regular vs Collet-Eckmann'' by Avila and Moreira in \cite{AM}.

In its foundation, the success with the Logistic Family relied on \emph{real a priori} bounds which are unavailable for genuine complex maps. Moreover, various \emph{near-Neutral} phenomena are not witnessed by the Logistic Family. For instance, there are Julia sets that are \emph{not} Locally Connected -- the \emph{non-JLC Category}. Even if the Julia set of a map is locally connected, the map can have many renormalization ``near non-JLC'' levels where ql-bounds fail.

Interest in the near-Neutral/Rotational Dynamics is rooted in the classical Celestial Mechanics. The first breakthroughs in the linearization of circle diffeomorphisms and local maps were made by Denjoy (1932) and Siegel (1942). In the 1950s and 60s, the Kolmogorov-Arnold-Moser (KAM) theory emerged, resulting in the reenvision of the near-rotation phenomenon in Mathematics and Physics. The Renormalization Ideas introduced in the mid-1970s led, in particular, to numerous conjectures in Low-Dimensional Dynamics, including the nature of the KAM small divisor problem. The theory of analytic circle diffeomorphisms and \emph{Local Theory} for neutral holomorphic germs and quadratic polynomials received an essentially complete treatment in the
second half of the last century in the work by Arnold (in the KAM framework), Cherry, Bruno,
Herman,  Yoccoz, and Perez-Marco.

For quadratic polynomials, the \emph{Neutral Family} 
\begin{equation}
\label{eq:neutr famil} z\mapsto e^{2\pi i \theta} z+ z^2\ \colon\ \C\ \selfmap, \qquad \theta\in \R/\Z,
\end{equation} 
encodes the boundary of the main hyperbolic component of $\MM$. There are two renormalization theories that analyze the \emph{Global Theory} of the Neutral Family.

\begin{figure}[t]  \centering
   \includegraphics[width=14cm,trim={0 5 0 5}, clip]{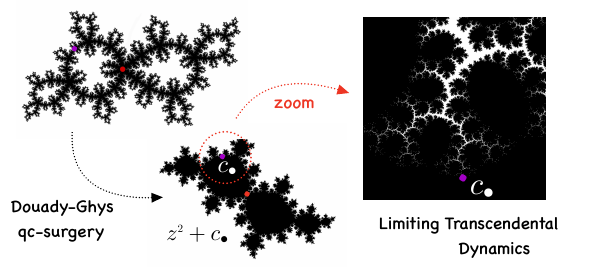}
   \caption{The 
   universality of the Golden Siegel map $z^2+c_\bullet$, where the parameter $c_\bullet$ is depicted in Figure~\ref{Fig:SiegelFamil}. The source of \emph{a priori} bounds for the respective renormalization theory is the  Douady-Ghys qc surgery; it transfers \emph{real a priori} bounds from Blaschke maps to complex Siegel polynomials of bounded-type.}
   \label{Fig:SiegelMaps}
\end{figure}

Renormalization Theory of (bounded-type) Siegel maps, also initiated by physicists, was mathematically designed by McMullen in~\cite{McMSiegel} in the mid-1990s. He established the dynamical universality as shown in Figure~\ref{Fig:SiegelMaps}. The parameter counterpart, Figure~\ref{Fig:SiegelFamil}, was developed two decades later in~\cite{DLS,DL_GAFA}, using the framework of the pacman renormalization. Theorem~\ref{thm:hyperb:pacman} states the hyperbolicity of the associated renormalization; its proof employs Transcendental Dynamics~\S\ref{ss:NearSiegel}.

\begin{figure}[t]  \centering
   \includegraphics[width=14cm,trim={0 8 0 8}, clip]{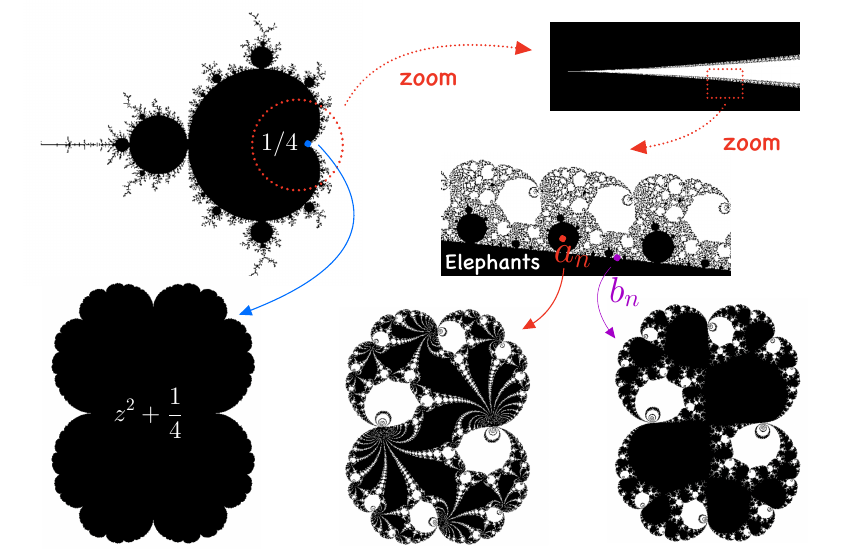}
   \caption{Parabolic Enrichment. Top: the Elephant structure in the parameter plane. Bottom: dynamical planes $p_{1/4}, p_{a_n}, p_{b_n}$. As $n\to \infty$, the pair $\big(p_{a_n}, p^{{\qq}_{1,n}}_{a_n}\big)$ converges to the commuting pair $\big(p_{1/4}, \ p^{[1]}_{a_\infty}\big)$, see Figure~\ref{Fig:ParabEnrich}, where $p^{[1]}_{a_\infty}$ is the associated Lavaurs transcendental map. Similarly, the neutral pairs $ (p_{b_n}, p_{b_n}^{{\qq}_{1,n}})$ converge to $(p_{1/4},\ p_{b_\infty}^{[1]})$. Lavaurs limits in the neutral case are efficiently encoded by continued fractions, see \S\ref{ss:SectorRenorm}.}
   \label{Fig:Cusp}
\end{figure}

The significance of the near-parabolic renormalization (see Figures~\ref{Fig:Cusp} and~\ref{Fig:ParabEnrich}) has been demonstrated by Shishikura in the 1990s in~\cite{Shi:ann}, where he showed that the Hausdorff dimension of $\partial \MM$ is equal to $2$. (The area of $\partial \MM$ is expected to be $0$, see~\S\ref{ss:MolecRenorm}.) The argument in \cite{Shi:ann} involved two iterates of the near-parabolic renormalization. Control of \emph{infinite iterates} became available in the mid-2000s when Inou and Shishikura established the hyperbolicity of the near-parabolic renormalization. New \emph{sectorial a priori} bounds, see Figure~\ref{Fig:SectorRenorm}, were put forward and justified in the perturbative regime; these were the first bounds in the \emph{non-JLC Category}. Many applications followed. For instance, the first examples of positive area Julia sets were constructed by Buff and  Ch{\'e}ritat in \cite{BC}, and then, of a different type, by Avila and Lyubich~\cite{AL}. In the mid-2010s, exact models of the postcritical sets were produced for neutral maps (Siegel and Cremer)~\cite{SY,Ch} in the Inou-Shishikura class.

Many global properties of the Neutral Family~\eqref{eq:neutr famil} are closely linked to those of unicritical circle maps; cf. the Douady-Ghys qc-surgery in Figure~\ref{Fig:SiegelMaps}. The Renormalization Theory of unicritical circle maps was fully developed (up to uniform hyperbolicity) in the works of de Faria, de Melo, and Yampolsky by the early 2000s.

In the mid-2000s, motivated by W. Thurston's work on $3$-manifolds and postcritically finite rational maps~\cite{BillTh, DH:Thurst}, Kahn introduced the Near-Degenerate regime to Renormalization theory~\cite{Kahn}. Jointly with Lyubich, they developed a machinery~\cite{KL09,KL2,KL3} that established quadratic-like \emph{a priori} bounds for ql-combinatorics of positive core-entropy $\ge \varepsilon >0$. It became apparent that this new theory provides an adequate replacement for the framework of real \emph{a priori} bounds in the genuine complex setting. However, at that time, the theory could not handle the zero-entropy setting (e.g., the non-JLC Category), and its development halted for more than a decade.

The Near-Degenerate regime in the zero-entropy setting was designed in~\cite{DL_Uniform}, see Theorem~\ref{thm:wZ}, using the framework of \emph{almost-invariant} pseudo-Siegel disks. Theorem~\ref{thm:wZ} established \emph{a priori} bounds for all maps in the Neutral Family~\eqref{eq:neutr famil}; in particular, providing the first \emph{non-perturbative a priori} bounds in the non-JLC Category. Various ideas behind~\cite{DL_Uniform} shaped follow-up research activities. In particular, Theorem~\ref{thm:MLC:Feigenb}, \cite{DLFeig} has finally confirmed the MLC at the classical period-doubling Feigenbaum parameter shown in Figures~\ref{Fig:Mand_set} and \ref{Fig:Feigenb:C}.

Advances in the renormalization theory of the Mandelbrot set have often provided guidance for more general settings. In the 2000s, various deep results for the Logistic Family have been extended to more general classes of real analytic maps; see \cite{ALdM, KSS,KSS-density}. For higher-degree polynomials, the Fibonacci-type renormalization discovered in~\cite{LM} does not degenerate; employing the real bounds, this renormalization was successfully understood, see~\cite{Sm2, Sm3}. With the introduction of the Near-Degenerate Regime, the complex Fibonacci case also became amenable to control, leading to generalizations of Yoccoz-type results to higher-degree polynomials~\cite{KLUnicr, AKLS}. Recently, in my joint work with Drach~\cite{DD}, Yoccoz-type realization and rigidity results were established for Transcendental Dynamics (for maps of controlled finite order and without asymptotic values) by developing a transcendental version of the near-degenerate analysis.

In the 2020s, the Near-Degenerate regime has also been introduced for rational maps. In my joint work with Luo~\cite{DLu}, \emph{a priori} bounds were established for hyperbolic components of disjoint type in the space of rational maps, leading to various realization results. This confirms McMullen's conjecture in the disjoint-type case from the 1990s, which asserts that Sierpi\'nski hyperbolic components are bounded. Simultaneously, Kahn~\cite{Kahn2} presented an independent approach with a different scope to \emph{a priori} bounds for establishing McMullen's conjecture (for all Sierpi\'nski cases). Meanwhile, Lim~\cite{Lim1} established bounds for rotating Herman curves of bounded-type rotation speed, leading to the development of the corresponding renormalization theory.

A popular version of the MLC story has been presented in a 2024 Quanta Magazine article \cite{Quest}.


\section{Renormalization operators.} \label{s:RrnormTheory:MM} 
\begin{figure}[t]  \centering
   \includegraphics[width=14cm,trim={0 10 0 8}, clip]{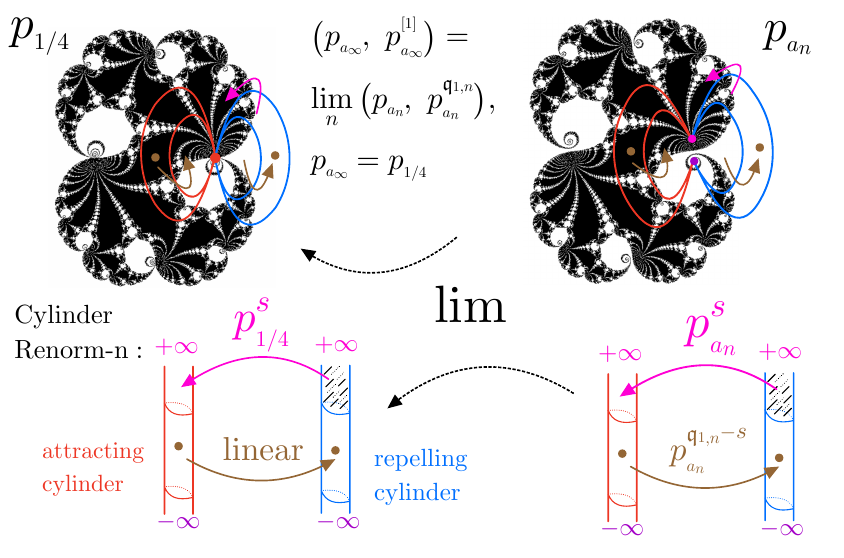}
   \caption{Top: as $a_n\to 1/4$, see Figure~\ref{Fig:Cusp}, the first return map $p^{\qq_{1,n}}_n$ converges to the transcendental Lavaurs map $p^{[1]}_{a_\infty}$ commuting with $p_{1/4}$. Bottom: the associated cylinder renormalization. See~\S\ref{ss:renorm:parablic} for details.}
   \label{Fig:ParabEnrich}
\end{figure}

Following the discovery of the period-doubling universality (Figure~\ref{Fig:Mand_set}), the framework of analytic renormalization operators acting on infinite-dimensional spaces of maps was put forward. The hyperbolicity of such operators always yields important consequences. There are two kinds of hyperbolicity frameworks:


\begin{itemize}
  \item[\setword{(QL)}{item:R:QL}] \emph{Quadratic-Like} Renormalization acting on ql maps as illustrated in Figure~\ref{Fig:QL_renorm}, top; and
  \item[\setword{(Neut)}{item:R:Cyl}] \emph{Near-Neutral} (Cylinder or Sectorial) Renormalization, see Figures~\ref{Fig:SectorRenorm} and~\ref{Fig:ParabEnrich}.
\end{itemize}

In addition to the~\ref{item:R:QL} and~\ref{item:R:Cyl} Cases, let us also recognize:
\begin{itemize}
      \item[\setword{(Puz)}{item:R:puz}] \emph{Puzzle} (or generalized ql) Renormalization encoding the first return map to a puzzle, Figure~\ref{Fig:QL_renorm}, bottom.
\end{itemize}

The~\ref{item:R:QL} bounded-type case is classical and now complete. The \emph{bounded-type} condition refers to parameters $c$ within nested intersection of small copies $\bigcap_n \MM_n$ as in~\eqref{eq:Yoccoz:results}, so that the relative periods of $\MM_{n+1}$ rel $\MM_n$ are uniformly bounded. Equivalently, in the dynamical planes of $p_c$, for every dynamical scale $\eta$ (as in~\S\ref{ss:resc:dynam}), $p_c$ has a quadratic-like structure $p^{\qq}_c\colon U'\to V'$ with a non-escaping critical orbit so that $|V'|\asymp \eta$. Similarly, the name ``bounded-type'' is used elsewhere in the renormalization theory, e.g., for neutral maps. In the late 1990s, the hyperbolicity problem in the bounded-type ql case was reduced to the justification of quadratic-like bounds. They have been recently supplied by Theorem~\ref{thm:MLC:Feigenb}, completing the program; see Theorem~\ref{thm:QL:hyperb}.

The hyperbolicity theory for neutral maps emerged in the mid-2000s, when Inou and Shishikura established the hyperbolicity for near-parabolic renormalization. We will state this result in~\S\ref{ss:NearParab}. The hyperbolicity of near-Siegel renormalization was established a decade later; see~\S\ref{ss:NearSiegel}. The neutral and satellite ql renormalization naturally form the \emph{Molecule renormalization}, see Figure~\ref{fig:mol map} and~\S\ref{ss:MolecRenorm}. The hyperbolicity of such a renormalization is likely required for Conjecture~\ref{cnj:hyp_comp} and the area-zero conjecture on $\partial \MM$. Conjecturally, there is no renorm-expanding regime on the Molecule: every parameter there has substantial hyperbolic components in all scales.

\begin{figure}[t!]  \centering
   \includegraphics[width=14cm, trim={0 10 0 8}, clip]{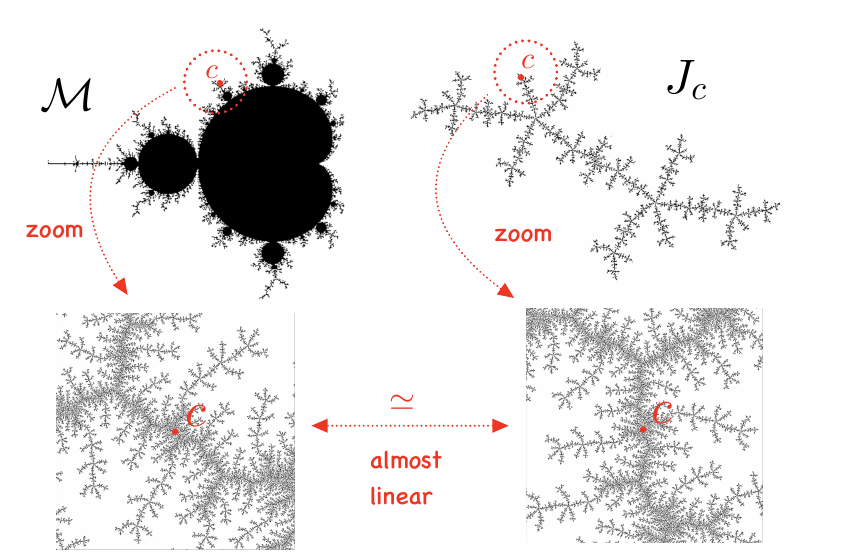}
   \caption{It is natural to ask if a form of Tan-Lei type similarity~\cite{Tan} (e.g., for the vein structure) between respective small neighborhoods of $\MM$ and $\Jul_c$ at $c$ occurs if and only if (parameter) hyperbolic components of $\MM$ become negligible; this is related to the notion of the renorm-expanding in~\S\ref{ss:resc:dynam}. For puzzles and parapuzzles, the \emph{Almost-Linear} relation (similarity) was established for Yoccoz maps in~\cite{L:Parap}, see~\S\ref{ss:puzzleRegime}.}
   \label{Fig:Exp_Case}
\end{figure}

In short, the Puzzle Renormalization deals with the renorm-expanding (see~\S\ref{ss:resc:dynam}) and ``virtual primitive ql'' (central cascades~\S\ref{ss:puzzleRegime}) regimes.  In the renorm-expanding regime, there are no geometric limits as in~\eqref{eq:resc:limits}, but there is a strong phase-parameter relation between the Mandelbrot and Julia sets as illustrated in Figure~\ref{Fig:Exp_Case}. One may hope that this case can also be neatly integrated into the hyperbolicity framework:

\begin{quest}[Uniform hyperbolicity rel $\partial \MM$]\label{conj:hyperb:main}
Is there a renormalization operator of a hyperbolic nature, associated with all parameters and scales of $\partial \MM$, with dimension $1$ unstable directions, such that the operator combines all the above kinds of renormalization regimes: ~\emph{\ref{item:R:QL},~\ref{item:R:Cyl},~\ref{item:R:puz}?}
\end{quest}
\noindent  Here is a hoped-for \emph{universality statement} from a desired operator in Question~\ref{conj:hyperb:main}. Consider a continuous zooming-in at $c\in\partial \MM$. If in a small $\delta$-scale of $c$, all hyperbolic components of $\MM$ are small rel $\delta$, then we expect a form of Tan-Lei type similarity between $\MM$ and $\Jul_c$ as depicted in Figure~\ref{Fig:Exp_Case}. If there is an unbounded (of size $\gg \delta$) hyperbolic component of $\MM$ hitting the $\delta$-scale of $c$, cf. Figure~\ref{Fig:SiegelFamil}, then we anticipate strong near-Neutral effects. And if all non-negligible hyperbolic components of $\MM$ are intermediate in the $\delta$-scale of $c$, cf. Figure~\ref{Fig:Feigenb:C}, then the quadratic-like features emerge. Finally, if the $\delta_1$ scale of $c_1$ is combinatorially close to the $\delta_2$-scale of $c_2$, then these scales should also be close geometrically. The last statement is a standard application of the hyperbolicity framework of the renormalization.

\subsection{Analytic and laminar structures on the space of ql maps.}\label{ss:RR:QL:analyt} 
\begin{figure}[t]  \centering
   \includegraphics[width=14cm,trim={0 0 0 0}, clip]{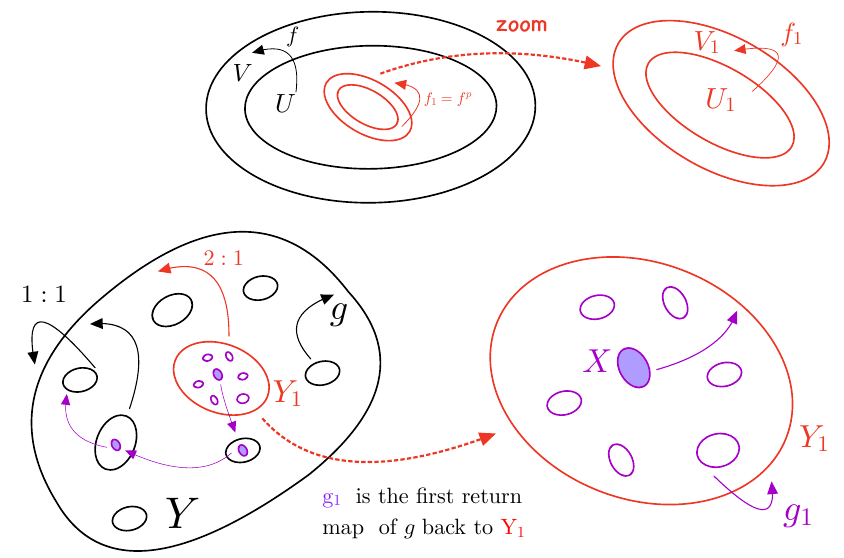}
   \caption{Top: Quadratic-Like Renormalization $f_1=f^p\colon U_1\ \overset{2:1}\longrightarrow\ V_1$ is obtained from $f\colon U \ \overset{2:1}\longrightarrow\  V$ by iterating, restricting, and then linear rescaling, see~\S\ref{ss:RR:QL:operator}. Bottom: Puzzle (or generalized ql) Renormalization is the first return map to a puzzle $Y_1$, see~\S\ref{ss:puzzleRegime}.}
   \label{Fig:QL_renorm}
\end{figure}

A \emph{quadratic-like map} $f\colon U \ \overset{2:1}\longrightarrow\ V$ is a degree $2$ holomorphic branched covering between open disks $U\Subset V\Subset \C$, see Figure~\ref{Fig:QL_renorm}. The set of its non-escaping points $ \filled_f\coloneqq \{z\in U\ : \ f^n(z)\in U \quad \forall n\ge1\}$
is called the \emph{filled Julia set} of $f\colon U\to V$. The boundary $\Jul_f\coloneqq \partial \filled_f$ is the \emph{Julia set} of $f$. The Julia set is connected if and only if the unique critical point of $f$ is non-escaping. In this case, $f$ is hybrid equivalent to a unique polynomial $p_c\colon z\mapsto z^2+c$ with $c=c(f)\in \MM$: there is a qc conjugacy $h$ between $f\colon U\to V$ and a quadratic-like restriction $p_c\colon A\to B$ such that $h$ is conformal on the non-escaping sets.

Here and elsewhere, by a slight thinning of $V$ and $U,$ we can assume that $\overline U$ is a closed topological disk. A \emph{small Banach neighborhood} of $f\colon U\to V$ is the set of all analytic maps close to $f$ in the $\sup$-norm:
\begin{equation}
\label{eq:dfn:BanachBall}
\mathcal B_\varepsilon (f\mid U) = \{g\colon U\to V \ :\ \|f-g\|_{\overline U}\ \le\  \varepsilon\}.    
\end{equation}
Let $V'\Subset V$ be a slight thinning of $V$. If $\varepsilon$ is sufficiently small, then every $g\in \mathcal B_\varepsilon (f\mid U)$ admits a quadratic-like structure $g\colon U_g\ \overset{2:1}\longrightarrow\ V'$, where $U_g$ is the unique $g$-lift of $V'$. This naturally introduces complex analytic structure for quadratic-like maps. There are two natural equivalences on the space of quadratic-like maps $f\colon U\to V$: linear conjugacy (i.e., rescaling) and thickening/thinning $U$ and $V$. In~\cite{LyuFeigen}, respecting these equivalences, the space of \emph{all} quadratic-like maps $\QL$ was supplied with complex-analytic structure induced by small Banach balls as in~\eqref{eq:dfn:BanachBall}. In this space, hybrid-equivalent maps form a codimension $1$ lamination with analytic leaves. The holonomy along hybrid leaves is qc and is asymptotically conformal. Such a holonomy combined with the hyperbolicity rigorously justifies various universalities as in Figure~\ref{Fig:Feigenb:C}. For near-Neutral Renormalization, establishing an appropriate version of the following consequence of the laminar structure is one of the key technical challenges, see~\S\ref{ss:NearSiegel}, \S\ref{ss:MolecRenorm}.

\begin{lem}[\cite{LyuFeigen}, $\codim=1$ stability of a ql $f$]\label{lem:LamStr} Let $f\colon U\to V$ be a ql map with a connected Julia set as above. Then, for a small $\varepsilon>0$, there is an analytic function $\tau\colon \mathcal B_\varepsilon (f\mid U) \to \C$ such that $f, g\in B_\varepsilon (f\mid U)$ are hybrid conjugate if and only if $\tau(f)=\tau(g)$. 
\end{lem}

\subsection{QL Compactness.}\label{ss:RR:QL:compact} Given a ql map  $f\colon U \ \overset{2:1}\longrightarrow\ V$, its \emph{modulus} is either 
\begin{equation}
\label{eq:bmod_f}
\mod (f)\coloneqq \mod( V\setminus U)\qquad \qquad\text{ or }\qquad\qquad \bmod(f)\coloneqq \mod (V\setminus \filled_f);
\end{equation}
both quantities are compatible. If $K_f$ is a Cantor set, then $\bmod(f)$ is the supremum of moduli over all annuli in the correct homotopy class separating $\partial V$ from $\filled_f$. Given a $\delta>0$, the space of ql maps with $\bmod(f)\ge \delta$ is compact: every sequence of ql maps $f_n$ with $\bmod(f_n)\ge \delta$ has a convergent subsequence within $\QL$; the limit $f_\infty$ also satisfies $\bmod(f_\infty)\ge \delta$. More generally, instead of separating $\partial V$ and $\filled_f$, we can measure the separation between $\partial V$ and the postcritical set. (In this framework, sequential limits of ql maps are, in general, outside the $\QL$ space.) This approach is adopted for near-neutral maps where the classical ql bounds fail, see~\S\ref{ss:SectorRenorm}.


\subsection{QL Renormalization operator.}\label{ss:RR:QL:operator} As Figure~\ref{Fig:QL_renorm} illustrates, a ql renormalization $f_1=f^p\colon U_1\ \overset{2:1}\longrightarrow\ V_1$ of $f\colon U\ \overset{2:1}\longrightarrow\ V$ is obtained by iterating $f$ followed by a ql-restriction and a linear rescaling. Since all three operations are analytic, the renormalization admits a natural analytic extension to nearby maps; we refer to such an extension as a \emph{renormalization operator}. Below, we first introduce a non-dynamical analytic operator. By design, it is always compact: the image of a ball is precompact. If there are appropriate geometric bounds, then the operator becomes  dynamical and hyperbolic.

 Let $U'\Subset U$ and $V'\subset V$ be slight thinning of $U$ and $V$ with $\delta\coloneqq \mod(U\setminus U')$. It is a fundamental fact (follows from Cauchy's integral formula) that the restriction operator 
\begin{equation}
\label{eq:compact:operator}
\mathcal B_\varepsilon (f\mid U)\ \to \  \mathcal B_\varepsilon (f\mid U'),\qquad\qquad [g\colon U\to \C]\ \mapsto \ [g\colon U'\to \C]
\end{equation}
is compact. Here $\varepsilon$ in $\mathcal B_\varepsilon (f\mid U)$, see~\eqref{eq:dfn:BanachBall}, is chosen sufficiently small depending on $\delta$. Viewing $f_1\equiv f^p\ \colon U_1\ \overset{2:1}\longrightarrow\ V_1$ as a subsystem of $f\colon U\ \overset{2:1}\longrightarrow\ V$,  we assume the following the \emph{compactness condition}:
\begin{equation}
\label{eq:dfn:delta}   \bigcup_{i=0}^{p-1} f^i(U_1)\Subset U'\qquad\qquad\text{ where }\qquad   \delta\coloneqq \mod(U\setminus U') \qquad \text{as above};
\end{equation}
i.e., the forward orbits of points in $U_1$ stay in $U'$ until their first return to $V_1\subset V'$. Then, since the restriction operator~\eqref{eq:compact:operator} is compact, we obtain a compact non-dynamical analytic operator 
\begin{equation}
    \label{eq:RR:QL}
\RR \colon \mathcal B_\varepsilon (f_1\mid U)\ \to  \  \mathcal B_{\varepsilon_2} (f^p\mid U_1),\qquad [g| U] \ \mapsto \ [g| U'] \ \mapsto \ [g_1=g^p|U_1].
\end{equation}

Suppose that $f\equiv f_0$ is infinitely ql-renormalizable: there is an infinite sequence $f_n\colon U_n  \overset{2:1}\longrightarrow\ V_n$ with $f_{n+1}=f_n^{p_{n+1}}$ as before. We also assume that $\delta>0$ is selected uniformly over $f_n$ and that the periods $p_n\le \overline p$ are uniformly bounded. If the sequence $f_n\colon U_n\to V_n$ has geometric bounds $\bmod (f_n)\ge \eta>0$, see~\eqref{eq:bmod_f}, then, by the compactness of ql maps, the $\omega$-limit set $\Hors_f$ of $(f_n)_n$ is compact. Such a limit set $\Hors$ is called a \emph{renormalization horseshoe}. If the geometric bounds have an additional ``beau'' property (it is often the case), then $\Hors_f$ is a hyperbolic set of $\RR$. Similarly, a renormalization horseshoe $\Hors_{\mathcal X}$ is defined for a set of maps $\mathcal X=\{f_\tau\equiv f_{0,\tau}\}$. It is natural to construct a horseshoe as big as possible. Applying Theorem~\ref{thm:MLC:Feigenb} together with~\cite{LyuFeigen}, we obtain:

\begin{thm}[Hyperbolicity, {\ref{item:R:QL}} Case]
\label{thm:QL:hyperb} For every bound $\overline p$ on renormalization periods, there is a hyperbolic renormalization horseshoe $\Hors_{\overline p}$ with $\dim=1$ unstable manifolds for the ql renormalization operator as above. The horseshoe $\Hors_{\overline p}$ attracts exponentially fast all infinitely renormalizable ql maps $f\equiv f_0$ with relative periods $p_n\le \overline p$.
\end{thm}

\begin{figure}[t]  \centering
   \includegraphics[width=13cm,trim={0 20 0 10},clip]{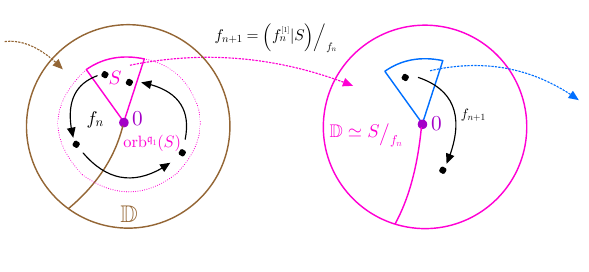}
   \caption{Sector Renormalization~\S\ref{ss:SectorRenorm} is the first return map to a sector $S$ whose sides are glued dynamically.}
   \label{Fig:SectorRenorm}
\end{figure}

\subsection{Sector Renormalization.}\label{ss:SectorRenorm} The Sector Renormalization was originally designed to study Local Dynamics around an indifferent fixed point by Douady and Yoccoz, with various ideas going back to the theory of Ecalle-Voronin’s invariants. For high-type rotation numbers, Inou and Shishikura extended the sector renormalization to a semi-global compact operator with geometric control of the critical orbit. Here, we provide more details on the semi-global case.

The idea of sector renormalization $\RR_{\sec}$ is illustrated in Figure~\ref{Fig:SectorRenorm}. In short, it is the first return map of $f_n$ to the appropriately chosen sector $S$ so that the critical orbit returns to $S$. Let us start with notation: 
 \begin{itemize}
 \item $\Theta=\{\theta=[0;a_1,a_2,\dots] \mid a_i \in \N_{\ge 1}\}\simeq  [0,1]\setminus \Q\simeq (\R/\Z)\setminus \Q$ is the set of irrational rotation numbers written as continued fractions;

  \item $\Theta_\bnd=\{\theta=[0;a_1,a_2,\dots] \mid a_i \le M_\theta\}\subset \Theta$ is the set of bounded (irrational) rotation numbers; 
  
  
 \item $\ovlTheta =\{\theta= [0;a_1,a_2,\dots] \mid a_i \in \N_{\ge 1}\cup \{\infty\}\}\supset \Theta$ is the formal compactification of $\Theta$; 
  \item $\Theta_{>N}=\{\theta=[0;a_1,a_2,\dots] \mid a_i \in \N_{> N}\}$; 
  \item $f_\theta \colon z\mapsto e^{2\pi i \theta}z +z^2$, $\theta\in \Theta$ is an irrationally-neutral quadratic polynomial as in~\eqref{eq:neutr famil},
   \item  $\pp_n/\qq_{n}\approx \theta$ is the best approximants of $\theta\in \Theta$ starting with $\qq_0=1$, 
   \item $f_{n}\equiv f_{n,\theta}\coloneqq \RR_{\sec}^n  f_\theta$ is the $n$th sector renormalization;
    \item if $\theta\in \Theta_{>1}$, then $\theta_n\coloneqq (-1)^n [0;a_{n+1},a_{n+2},\dots]$ is the rotation number of $f_n$ and\\  $f_\theta^{[n]}\coloneqq f^{\qq_n}_\theta$ is the sector prerenormalization of $f_n= f^{[n]}_\theta/_{f^{[n-1]}_\theta}$. 
 \end{itemize}
The definitions of $\theta_n$ and $f^{[n]}_\theta$ require a minor adjustment for the general case ``$\theta\in \Theta$''. We view $f_n=\RR^n_{\sec} f_\theta$ as 
\begin{equation}
    \label{eq:dfn:f_n}f_n\colon \Disk\dashrightarrow \Disk\qquad\qquad f_n(0)=0, \qquad f_n(0)=e^{2\pi i \theta_n}, \qquad \theta_n\in\Theta.
\end{equation}
The \emph{sector renormalization} $f_{n+1}$ of $f_n$ as shown in Figure~\ref{Fig:SectorRenorm} is the first return map to the smallest sector $S$ bounded by curves $\gamma$ and $f_n(\gamma)$ emerging from $0$ so that the sides $\gamma,  f_n(\gamma)\subset S$ are dynamically glued by $f_n$: 
\begin{equation}
    \label{eq:sector:gluing} \gamma \ni x \ \sim \ f_n(x) \in f_n(\gamma)\qquad\qquad \text{ so that }\qquad S/_{f_n} \simeq \Disk.
\end{equation}
 We assume that $\gamma$ lands at $0$ at a well-defined angle. The ``smallest'' condition on $S$ refers to the requirement that $S$ contains the smallest angle between $\gamma$ and $f_n(\gamma)$. The sector $S$ should be chosen appropriately so that the first return map to $S$ controls the critical orbit. Before taking the quotient~\eqref{eq:sector:gluing}, the first return map is a commuting pair of the form $F_n=\big(f_n^{[1]}\equiv f^{\qq_1(\theta_n)}_n,\ f^{\qq_1(\theta_n)\pm 1}_n\big).$ The quotient map by $f_n$ identifies both branches in the commuting pair $F_n$; we can write $f_{n+1}= \RR_\sec f_n=f_n^{[1]}/_{f_n}\equiv F_n/_{f_{n+1}}$, and we refer to both $F_n$ and $f_n^{[1]}$ as the \emph{sector prerenormalizations} of $f_{n+1}$. Similarly, $f^{[n]}_\theta$ is the $n$th sector prerenormalization of $f_n$, and we write: $f_{n}= f^{[n]}_\theta/_{f^{[n-1]}_\theta}$. As in the quadratic-like case, see~\eqref{eq:dfn:delta}, we require a \emph{compactness condition}:
\begin{equation}
    \label{eq:dfn:delta:sect} \mod\big(\Disk\setminus \orb^{\qq_1} (S)\big)\ge \delta, 
\end{equation}
where $\orb^{\qq_1} (S)$ is formed by orbit trajectories of points in $S$ projecting to $\operatorname{dom} (f_{n+1})$. We refer to~\eqref{eq:dfn:delta:sect} as \emph{sectorial bounds}; they imply the compactness of $\RR_\sec$. As an application of Theorem~\ref{thm:wZ}, the compact sectorial renormalization operator $\RR_\sec$ is justified in~\cite{DL:sector-bounds} for all neutral quadratic polynomials $(f_\theta)_{\theta\in \Theta}$ with a compactness condition~\eqref{eq:dfn:delta:sect} for a uniform $\delta$ and with well-defined parabolic limits parametrized by $\ovlTheta$.

\subsection{Parabolic Enrichment and the Cylinder Renormalization.}\label{ss:renorm:parablic} We recall that the filled Julia set is discontinuous at parabolic parameters: as shown in Figure~\ref{Fig:Cusp}, the limit of filled Julia sets can be smaller than the filled Julia set of the limiting parameter. This discontinuity can be resolved by introducing limiting Lavaurs maps. Let us illustrate a general idea with examples.

Consider a sequence of parameters $a_n\to 1/4$ as shown in Figure~\ref{Fig:Cusp}. In the dynamical plane of $p_{a_n}\colon z\mapsto z^2+a_n$, see Figure~\ref{Fig:ParabEnrich}, the critical orbit returns after $\qq_{1,n}\equiv \qq_1(a_n)$ iterates. The orbit spends $\qq_{1,n}-s$ iterations traveling from the red cylinder to the blue cylinder, and then $s$ iterations traveling back. We have a commuting pair $(p_{a_n}, p_{a_n}^{\qq_{1,n}})$. In the limit, the parabolic gap is closed, and  $(p_{a_n},\  p_{a_n}^{\qq_{1,n}})$ converges to the pair $(p_{1/4}=p_{a_\infty},\ p_{a_\infty}^{[1]})$, where $p_{a_\infty}^{[1]}$ is a transcendental \emph{Lavaurs} map defined on the attracting basin of the parabolic point of $p_{1/4}$. The red and blue cylinders converge to the limiting attracting and repelling cylinders of $p_{1/4}$.

The \emph{cylinder renormalization} $\RR_\cyl$ is a convenient and canonical extension of a sector renormalization $\RR_{\sec}$ from~\S\ref{ss:SectorRenorm} such that the renormalization sector of $\RR_{\sec}$ is \emph{extended to a cylinder} of $\RR_\cyl$. While $\RR_{\sec}$ is defined up to conformal conjugacy, $\RR_\cyl$ is defined up to a linear conjugacy. For $p_{a_n}$, its \emph{cylinder renormalization} $\RR_\cyl (p_{a_n})\simeq p_{a_n}^{\qq_{1,n}}/_{p_{a_n}}$ is the first return map of $p_{a_n}$ back to the red cylinder, where its boundaries are identified using $p_{a_n}$. The limiting cylinder renormalization $\RR_\cyl\big(p_{1/4},\ p_{a_\infty}^{[1]}\big)$, also known as \emph{parabolic}, is the quotient of $p_{a_\infty}^{[1]}$ under $p_{1/4}$. We write
 \begin{equation}
 \label{eq:Renorm:Parab}
 \RR_\cyl (p_{a_n}) \ \simeq \ p^{\qq_{1,n}}_{a_n}\big/_{p_{a_n}}   \ \ \longrightarrow \ \ \RR_\cyl\big(p_{1/4},\ p_{a_\infty}^{[1]}\big)\simeq p^{[1]}_{1/4}\big/_{p_{1/4}},
 \end{equation}
 and refer to $\RR_\cyl (p_{a_n})$ as \emph{near-parabolic} for $n\gg 1$. In fact, $p_{a_\infty}^{[1]}$ and $\RR_\cyl\big(p_{a_\infty}^{[1]}\big)$ are uniquely determined by a linear map, called a \emph{transitional map}, between red and blue cylinders, see the bottom-left part of Figure~\ref{Fig:ParabEnrich}. Conversely, given a transitional map, one can construct the associated Lavaurs map and then realize it as a limit of nearby maps. A fundamental question is the dynamical control of this \emph{parabolic enrichment}.  As it was discussed in~\S\ref{ss:SectorRenorm}, renormalization of neutral quadratic polynomials is parametrized by $\theta= [0;a_1,a_2,\dots]\in \ovlTheta$. Parabolic renormalization applies to the cases where $a_i=\infty$. For illustration, see parameters $b_n$ in Figure~\ref{Fig:Cusp}.

\subsection{Hyperbolicity of the near-parabolic renormalizations.} \label{ss:NearParab} In~\cite{IS}, the hyperbolicity of the cylinder renormalization has been established for high-type parameters $\theta\in \Theta_{\le N},$ where $ N\gg 1$. The hyperbolicity implies the sectorial bounds as they were introduced in~\S\ref{ss:SectorRenorm}.  As we mentioned in Section~\ref{s:History}, this result had many applications. The strategy in \cite{IS} is as follows. First, with computer assistance, the \emph{a priori} bounds (a variant of~\eqref{eq:dfn:delta:sect}) are established for the parabolic renormalization representing 
\[  [0;\infty, \infty,\infty,\dots]\in \ovlTheta,\qquad \qquad \text{ see \S\ref{ss:SectorRenorm}.}\]
Using Teichm\"uller contraction, the geometric bounds are promoted to the hyperbolicity of the associated parabolic renormalization. By the stability of hyperbolicity, it also holds for all $\theta\in\Theta_{> N}$ for some $N\gg 1$.

\subsection{Hyperbolicity of near-Siegel Renormalization.} \label{ss:NearSiegel} Following notations of~\ref{ss:SectorRenorm}, consider maps $f_\theta \colon z\mapsto e^{2\pi i \theta}z +z^2$. For $\theta\in\Theta_\bnd$, the Douady-Ghys surgery transfers \emph{real bounds} from circle maps (the Blaschke family) to quadratic polynomials, see Figure~\ref{Fig:SiegelMaps}, left. Consequently, $f_\theta$ has a qc Siegel disk if and only if $\theta\in\Theta_\bnd$. Moreover, such a qc Siegel disk contains a critical point on the boundary. Relying on such input,  McMullen in~\cite{McMSiegel} developed dynamical renormalization theory for $\theta\in \Theta_\bnd$, and justified the dynamical universality: the exponential convergence of renormalizations towards the limit, see Figure~\ref{Fig:SiegelMaps}, right.

For two decades, it remained an open question whether dynamical universality could be promoted to hyperbolicity. The key obstacle was to establish a version of Lemma~\ref{lem:LamStr}: proving that Siegel maps remain stable under small perturbations, provided the perturbation does not change the multiplier of the Siegel fixed point. Here, perturbations are taken within Banach balls (cf.~\eqref{eq:dfn:BanachBall}) consisting of maps defined in neighborhoods of their Siegel disks. This obstacle was resolved by employing transcendental techniques on unstable manifolds. \emph{A priori} bounds imply the compactness of the operator, thus unstable manifolds exist but may have a large dimension. It was shown in~\cite{DLS} that maps on the unstable manifolds admit maximal analytic extensions as $\sigma$-proper maps ``$W\to\C$,'' cf.~\eqref{eq:resc:limits}, where $W$ is open and dense, with a single critical orbit. Such global $\sigma$-proper maps have $\codim=1$ stability as in Lemma~\ref{lem:LamStr}; this stability is then promoted to $\codim=1$ stability in the original class of maps, and subsequently to:

\begin{thm}[\cite{DLS}]
\label{thm:hyperb:pacman}
Pacman/Siegel Renormalization operator is hyperbolic with $\dim=1$ unstable manifold for parameters in $\Theta_\bnd$ of periodic type.
\end{thm}

\noindent For the Golden-mean Siegel map shown in Figure~\ref{Fig:SiegelMaps}, the hyperbolicity was established earlier by Gaidashev and Yampolsky with computer-assistance.  We also remark that the global $\sigma$-proper structure of renormalization fixed points (both of the quadratic-like and Siegel types) was discovered by McMullen in the 1990s. The theory of external transcendental rays for unstable manifolds was designed in~\cite{DL_GAFA}, leading to further applications.

\begin{figure}
  \begin{center}
   \includegraphics[width=14cm,trim={10 38 10 20},clip]{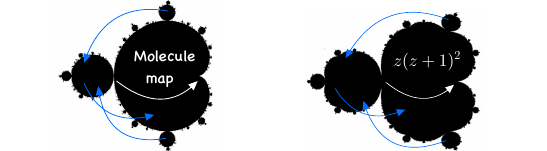}
 \end{center}
\caption{Left: the Main Molecule $\Mol$ of the Mandelbrot set and the induced by the Branner-Douady surgery Molecule map, see~\S\ref{ss:MolecRenorm}. The Molecule map is modeled by $z\mapsto z(z+1)^2$, right.}
\label{fig:mol map}
\end{figure}

\subsection{Molecule Renormalization.}\label{ss:MolecRenorm} As we have indicated in~\S\ref{ss:SectorRenorm}, Theorem~\ref{thm:wZ} implies the existence of a compact sectorial renormalization operator $\RR_\sec$ rel $\ovlTheta$. The proof of Theorem~\ref{thm:hyperb:pacman}, relying on the analysis of the transcendental dynamics on the unstable manifolds, provides a strategy towards:


\begin{conj}[Hyperbolicity of Neutral Renormalization]\label{conj:hyp:neutral} There is a hyperbolic cylinder renormalization operator $\RR_\cyl$ with $\dim=1$ unstable manifolds for all parameters in $\ovlTheta$.
\end{conj}

 We recall from~\cite[Appendix C]{DLS} that the satellite and neutral renormalizations are naturally unified into the Molecule Renormalization; its slowed-down version is illustrated in Figure~\ref{fig:mol map}. A natural extension of Conjecture~\ref{conj:hyp:neutral}:

\begin{conj}[Hyperbolicity of Molecule Renormalization]\label{conj:hyp:Molecule} There is a hyperbolic renormalization with $\dim=1$ unstable manifolds parametrized by the sped-up Molecule map as in Figure~\ref{fig:mol map}.
\end{conj}
The completion of Problem~\ref{prob:Satel:unbounded} is required for Conjecture~\ref{conj:hyp:Molecule}. As we have already mentioned, Conjecture~\ref{conj:hyp:Molecule} is likely required for the satellite cases of Conjecture~\ref{cnj:hyp_comp} and the Area Problem of $\partial \MM$.

\subsection{Puzzle Renormalization.} \label{ss:puzzleRegime} 
Puzzle techniques were first employed by Branner and Hubbard in~\cite{BrHub} to analyze cubic polynomials with one escaping critical point. In that setting, puzzles do not cut the Julia set. A few years later, Yoccoz introduced puzzles for quadratic polynomials and used them to reduce the MLC to the infinitely ql renormalizable case, see~\eqref{eq:Yoccoz:results}. Yoccoz puzzles cut a Julia set into smaller pieces. Then the Gr\"otzsch inequality is used to show that such pieces shrink to points.

In~\cite{LyuActa,L:Parap}, in the 1990s, the puzzle techniques were incorporated into the renormalization framework as shown in Figure~\ref{Fig:QL_renorm}. Given a puzzle $Y$ around a critical point so that the forward $p_c$-orbit of $\partial Y$ is disjoint from $Y$, we can consider the first return map $g$ of $p_c$ to $Y$. The domain of $g$ consists of countably many puzzles. One of them, say $Y_1$, contains a critical point and is mapped with degree $2$ into $Y$. All other pieces are mapped with degree $1$. The \emph{renormalization} $g_1$ of $g$ is the first return map of $g$ into $Y_1$. Let $X\Subset Y_1$ be the next-level puzzle containing the unique critical point. If the return time of $X$ rel $g$ is bigger than $1$, i.e., $X$ travels around $Y$ under $g$, then the ``asymmetric modulus'' grows; roughly, $g_1\colon Y_1\dashrightarrow Y_1$ has more expansion than $g\colon Y\dashrightarrow Y$. On the other hand, $g_1$ belongs to a \emph{central cascade} if the return time of $X$ is $1$, i.e., $g\mid X=g_1\mid X$.

If a central cascade is long but still finite, then the associated dynamical scales will be dominated by a certain primitive small copy $\MM'$ of $\MM$ that does not contain $c$. We refer to this phenomenon as \emph{virtual primitive ql}. 

If $g\colon Y\dashrightarrow Y$ is sufficiently expanding (i.e., moduli are big, cf. the renorm-expanding case in~\S\ref{ss:resc:dynam}), then there is an \emph{Almost-Linear puzzle-parapuzzle} relation~\cite{L:Parap}; i.e., the geometry of (puzzles of) $g\colon Y\dashrightarrow Y$ can be efficiently transferred to the parameter plane, cf. Figure~\ref{Fig:Exp_Case}. This is an essential ingredient for the ``regular vs stochastic'' dichotomy. Purely renorm-expanding real parameters (\S\ref{ss:resc:dynam}) are precisely Lebesgue density points for stochastic parameters: most of their small perturbations have exponential growth of the expansion along the critical orbit (the Collet-Eckmann condition), i.e., they are stochastic.

\section{Near-Degenerate Regime.}\label{s:NDR} As we mentioned in~\S\ref{s:History}, the Near-Degenerate Regime provides an adequate framework towards establishing \emph{a priori} bounds. Since it operates at a basic topological level, arguments developed for different dynamical scales are often compatible and can be merged provided
the setting is appropriately prepared. The non-dynamical foundations of the theory have been developed in~\cite{KL09}. Below, we will omit most of the technical discussion and indicate only the sources of geometric control.

\begin{figure}[t]  \centering
   \includegraphics[width=14cm,trim={0 15 0 10},clip]{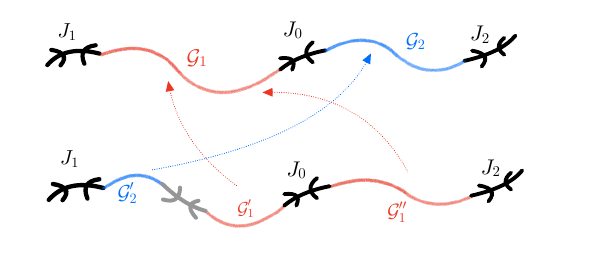}
   \caption{Illustration to the positive entropy argument, \S\ref{ss:NDR:PositiveEntr}, the airplane combinatorics.}
   \label{Fig:NDG:PositiveEntropy}
\end{figure}

On the dynamical scales with definite \emph{positive core entropy} $\ge \varepsilon>0$, the theory is fully developed, see \S\ref{ss:NDR:PositiveEntr} and \S\ref{ss:MLC:puzzle}. As Figure~\ref{Fig:NDG:PositiveEntropy} illustrates, the positivity of entropy prevents dynamically invariant horizontal degeneration. Thus, if there is a substantial degeneration, then it can be amplified at shallower levels.

\begin{figure}[t]  \centering
   \includegraphics[width=14cm,trim={0 15 0 10},clip]{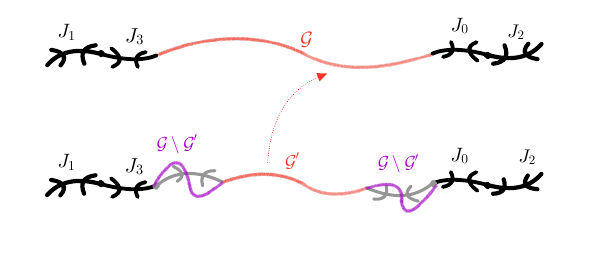}
   \caption{Illustration to~\S\ref{ss:MLC:QL:sat}: if $\mathcal G$ is invariant rel $\delta$, then $\Width(\mathcal G\setminus \mathcal G')\ge C_\delta\Width(\mathcal G)\gg \Width(\mathcal G)$.}
   \label{Fig:NDR:ZeroEntropy}
\end{figure}

Consider the \emph{zero core entropy} scales. On satellite ql levels, degeneration along non-invariant rectangles is amplified at shallower levels, while degeneration along combinatorially invariant rectangles is amplified at deeper levels; see Figure~\ref{Fig:NDR:ZeroEntropy} and~\S\ref{ss:MLC:QL:sat}. In the (near-)neutral case, invariant degeneration can actually exist within parabolic fjords because of the non-JLC; \emph{pseudo-Siegel} disks are constructed to consume most of the invariant degeneration, see the left loop of Figure~\ref{Fig:NearDegRegimIllustr}. Its right loop illustrates a strategy to develop pseudo-Siegel disks in the remaining near-neutral satellite cases. If achieved, this would fully unify the cases of \S\ref{ss:NDR:PositiveEntr}--\S\ref{ss:MLC:QL:neut} and confirm, in particular, the Satellite Case of the MLC.

Let us now discuss the \emph{Interpolation Problem} between the positive and zero entropy settings. It concerns parameters $c\in \MM$ that are geometrically close to the main Molecule $\Mol$ of $\MM$ but are outside of satellite little copies centered at $\Mol$. In the dynamical plane of $p_c\colon z^2+c$, there are many geometric scales of zero-entropy nature associated with $\Mol$. A key issue is that, on these scales, the postcritical set of $p_c$ can look as if it is dense in the Julia sets; i.e., the entropy on the postcritical set can blow up in the limit as $c$ approaches $\Mol$. (The postcritical set is only lower semicontinuous.) This effectively compels us to deal with \emph{virtual Julia sets} containing the critical orbit up to several first returns. We briefly formulate this \emph{Virtual Renormalization} setting in~\S\ref{ss:virtual_renorm}.

 \subsection{Primitive QL levels.}\label{ss:NDR:PositiveEntr} Assume that $f_1$ is a primitive ql renormalization of $f$ with a period $p$ as in Figure~\ref{Fig:QL_renorm}.  Instead of $\bmod(f)=\mod(V\setminus \filled_f)$ as in~\eqref{eq:bmod_f}, it is more convenient to consider the dual quantity, called the \emph{width} or the \emph{degeneration} of $f$:
\begin{equation}  
     \Width(f) \coloneq \frac{1}{\bmod(f)} \qquad\qquad  \text{so that }\qquad \bmod(f)\le \varepsilon \qquad\text{ iff }\qquad  \Width(f)\ge K\equiv \frac{1}\varepsilon,
\end{equation}
and similar for $f_1$. The width $\Width(f)$ is the extremal width of the family of curves connecting $K_f$ to $\partial V$.

Let us analyze how the degeneration of $f_1$ is developed in the dynamical plane of $f$. Since $f$ has $p$ small (filled) Julia sets of $f_1$, the total degeneration around these small Julia sets is $\Width^\tot\asymp p \Width(f_1)-O_p(1)$, where an additive error $O_p(1)$, that depends on $p$, comes from a variant of the Thin-Thick Decomposition. We can decompose $\Width^\tot=\Width^\ver+\Width^\hor$, where
\begin{itemize}
    \item  $\Width^\hor$ is the \emph{horizontal} degeneration between any pair of two small Julia sets; and
    \item  $\Width^\ver$ is the \emph{vertical} degeneration between $\partial V$ and any small Julia set.
\end{itemize}
 A key claim in~\cite{Kahn} is that $\Width^\ver\succeq \Width^\tot-O_p(1)$; i.e., the vertical part has a definite proportion. Since small Julia sets of $f_1$ are within the bigger Julia set of $f$, the claim implies
\begin{equation}
    \label{eq:key:PosEnropy} \Width(f)\ge \Width^\ver \ \succeq\  p\Width(f_1) -O_p(1).
\end{equation}

The idea behind $\Width^\ver\succeq \Width^\hor$ is illustrated in Figure~\ref{Fig:NDG:PositiveEntropy}. Assume by contradiction that  $\Width^\hor\gg \Width^\ver$, i.e., $\Width^\ver$ is negligible. It is a meta-principle that a degeneration should be invariant under dynamics. Therefore, the degeneration $\Width^\hor$ is (eventually) aligned with the Hubbard tree of $f$ and is represented by two \emph{wide rectangles} $\mathcal G_1$ and $\mathcal G_2$ depicted on the figure as narrow (for convenience) red and blue rectangles in the pencil style. Since $\Width^\ver$ is negligible, 
 \begin{itemize}
     \item (most of) $\mathcal G_2$ overflows a lift $\mathcal G''_1$ of $\mathcal G_1$; and 
     \item (most of) $\mathcal G_1$ overflows a lift $\mathcal G'_2$ of $\mathcal G_2$ and then a lift of  $\mathcal G'_1$ of $\mathcal G_1$.
 \end{itemize}
This is impossible by the  Gr\"otzsch inequality and~\eqref{eq:key:PosEnropy} follows.

By applying~\eqref{eq:key:PosEnropy} to $f$ and its $n$th renormalization $f_n$ with a big period $p_n\gg 1$; we obtain
\begin{equation}
    \label{eq:key:PosEnropy:f_n} \Width(f)\ \succeq \ p_n\Width(f_n) -O_{p_n}(1)\ \gg\  \Width(f_n)\qquad\qquad \text{if }\quad \Width(f_n)\gg_{p_n}\ 1.
\end{equation}
This implies \emph{a priori} bounds in the primitive bounded-type ql case: the degeneration $\Width(f_m)$ can not become big because, otherwise, the degeneration of $\Width(f_{m-n})$ will be even bigger. 

\subsection{Puzzle levels.}\label{ss:MLC:puzzle} We emphasize the following subtlety in~\eqref{eq:key:PosEnropy}: the period of $p_n$ is assumed to be big but \emph{still bounded} as the additive error $O_{p_n}(1)$ depending on $p_n$ needs to be bounded.

Assume now that the period $p$ of $f_1$ rel $f$ is big. We assume here that $f_1$ is the first ql renormalization of $f$ of the primitive type; this renormalization is encoded by a maximal primitive copy $\MM_1$ of $\MM$. The renormalization $f\mapsto f_1$ can be factorized using the puzzle renormalization, see Figure~\ref{Fig:QL_renorm}, 
\begin{equation}
\label{eq:PuzRenorm}
  f\mapsto g_1\mapsto g_2\mapsto \dots \mapsto g_k \mapsto f_1.
\end{equation}
If the copy $\MM_1$ is $\varepsilon$-away from the main molecule (the \emph{anti-Molecule} property, cf. Figure~\ref{fig:mol map}) and $p\gg_\varepsilon \ 1$, then the (modified) Principal Nest of~\eqref{eq:PuzRenorm} is combinatorially substantial and provides amplification of the degeneration. Roughly, the argument goes as follows. Using the notations of Figure~\ref{Fig:QL_renorm}, as $X$ travels around $Y$, the~\cite[Covering Lemma]{KL09} efficiently spreads the degeneration around. Then the~\cite[Quasi-Additivity Law]{KL09} collects the degeneration along the orbit of $X$; this eventually leads to:
\begin{equation}
    \label{eq:f:f_1:PrincNest} \Width(f) \gg \Width (f_1)\qquad\qquad \text{ if }\quad \Width(f)\gg_\varepsilon\ 1\quad \text{ and }\quad  p\gg_\varepsilon 1.
\end{equation}
Together,~\eqref{eq:key:PosEnropy:f_n} and~\eqref{eq:f:f_1:PrincNest} imply the \emph{a priori} bounds and the MLC under the anti-molecule condition~\cite{KL2,KL3}.

\begin{figure}[t]  \centering
   \includegraphics[width=16cm,trim={0 5 0 5},clip]{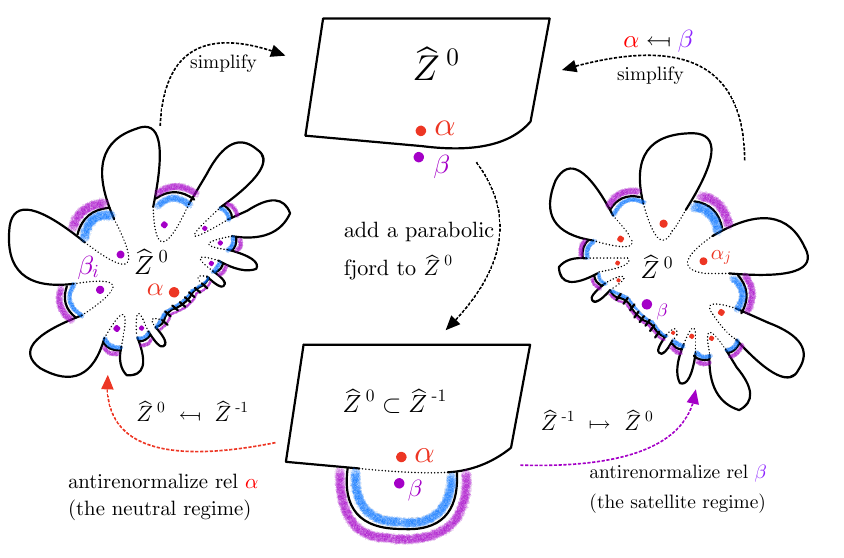}
   \caption{Left loop: the interaction between pseudo-Siegel disks and neutral sector antirenormalization, see~\S\ref{ss:MLC:QL:neut}. Right loop: illustration to the strategy to extend the theory of pseudo-Siegel disks to the satellite case.} 
   \label{Fig:NearDegRegimIllustr}
\end{figure}

\subsection{Satellite QL levels.}\label{ss:MLC:QL:sat} Assume now that $f_1$ is a satellite ql renormalization of $f$ of bounded-type. The key difference with the primitive case~(\S\ref{ss:NDR:PositiveEntr}) is the existence of combinatorially invariant rectangles $\mathcal G$ between small Julia sets as illustrated in Figure~\ref{Fig:NDR:ZeroEntropy}. Below, we sketch how to continue the argument behind~\eqref{eq:key:PosEnropy:f_n} in the presence of such a $\mathcal G$.

Let us assume that most of the $\Width^\tot$ is within $\mathcal G$. For a small $\delta>0$, the $(1-\delta)$ part of $\mathcal G$ overflows its lift $\mathcal G'$ as shown in the figure. However, before entering the rectangle $\mathcal G'$, curves forming the above $(1-\delta)$-part of $G$ must travel through preperiodic Julia sets as depicted using purple color. Let us denote by $\mathcal G\setminus \mathcal G'$ the family of such purple subcurves. Since $\mathcal G, \mathcal G'$ have almost the same width, up to ``$(1-\delta)$'', the Gr\"otzsch inequality implies
\begin{equation}
    \label{eq:NDR:SatCase} \Width(\mathcal G\setminus \mathcal G') \ge C_\delta\  \Width(\mathcal G), \qquad\qquad \text{where $\ C_\delta\to \infty\ $ as $\ \delta\to 0$}.
\end{equation}

The~\cite[Wave Lemma]{DLFeig} allows us to interpret the quantity $\Width(\mathcal G\setminus \mathcal G')$ as the degeneration  $\Width(f_{2})$ on a deeper level. Considering now $f, f_n$ and $f_{n+1}$, Estimate~\eqref{eq:NDR:SatCase} gives us the following refinement of 
\eqref{eq:key:PosEnropy:f_n}:  either 
\begin{enumerate}[label=\text{(\roman*)},font=\normalfont,leftmargin=*]
    \item\label{item:NDR:sat:1} $ \Width(f)\ \succeq_\delta \ p_n\Width(f_n) -O_{p_n}(1)$; or
    \item\label{item:NDR:sat:2} $\Width(f_{n+1})\ge C_\delta\ \Width(f_n)$.
\end{enumerate}

The Alternative~\ref{item:NDR:sat:2} is prevented by the Teichm\"uller contraction. Indeed, given a bound $\overline p$ on the relative period $p_{m+1}/p_m$ of $f_{m+1}$ rel $f_m$, the degeneration of $\Width(f_m)$ can grow at most exponentially fast:
$\Width(f_m)=O(\Delta^m_{\overline p})$. Therefore, by selecting $\delta>0$ to be sufficiently small so that $C_\delta\gg \Delta_{\overline p}$ is sufficiently big, Alternative~\ref{item:NDR:sat:2} is eventually prevented on some levels. We also select $n$ sufficiently big so that $p_n$ dominates the multiplicative constant representing ``$\succeq_\delta$''.  Alternative~\ref{item:NDR:sat:1} implies \[ \Width(f)\ \succeq_\delta \ p_n\Width(f_n) -O_{p_n}(1)\ \gg\  \Width(f_n)\qquad \qquad \text{ if }\  \Width(f_n)\gg_{p_n} \ 1.\]

\begin{thm}[\cite{DLFeig},\ \S\ref{ss:NDR:PositiveEntr}+\S\ref{ss:MLC:QL:sat}]
\label{thm:MLC:Feigenb} A priori bounds and the MLC hold for all ql-combinatorics of bounded-type. 
\end{thm}
\noindent In particular, the MLC holds at the period-doubling Feigenbaum parameter illustrated in Figure~\ref{Fig:Feigenb:C}. The arguments of~\S\ref{ss:MLC:puzzle} can also be integrated with those in Theorem~\ref{thm:MLC:Feigenb}; the combination, see~\cite{DKL}, yields the MLC along a real line and along any vein of the Mandelbrot set.

\subsection{Near-Neutral levels.} \label{ss:MLC:QL:neut} Let us first detail the following statement:
\begin{thm}[\cite{DL_Uniform}]
\label{thm:wZ}
Uniform pseudo-Siegel a priori bounds hold for all neutral quadratic polynomials.
\end{thm}

Using the notations of~\S\ref{ss:SectorRenorm}, consider a neutral polynomial $f_\theta \colon z\mapsto e^{2\pi i \theta}z +z^2$ for an irrational $\theta\in \Theta$. Theorem~\ref{thm:wZ} asserts that it has a sequence of nested pseudo-Siegel disks 
\begin{equation}
    \label{eq:dfn:wZ^m:intro}\qquad
\wZ^{-1}\supseteq \wZ^{0}\supseteq \wZ^{[1]}\supseteq \dots \supset \wZ^{[m]}\supseteq \dots \supseteq H_{f_\theta}\ \ni 0, \qquad \qquad \bigcap_{m\ge -1} \wZ^{[m]} =H\equiv H_{f_\theta},
\end{equation}
where $H$ is the \emph{Mother Hedgehog} (the filled postcritical set) of ${f_\theta}$ such that every $\wZ^{[m]}$ is almost invariant up to $\qq_{[m+1]}$-iterations. The disks $\wZ^{-1}$ and $\wZ^0$ are uniformly qc over all $f_\theta$ with $\theta\in \Theta$.

Under the $n$th sector renormalization $f_\theta \mapsto f_n$, cf.~Figure~\ref{Fig:SectorRenorm}, the pseudo-Siegel disk $\wZ^{[n+m]}_{f_\theta}$ of $f_\theta$ becomes $\wZ^{[m]}_{f_n}$ of $f_n$; i.e., every $\wZ^{[m]}$ becomes uniformly qc after applying the $m$th renormalization.

Figure~\ref{Fig:NearDegRegimIllustr}, the left loop, illustrates the interaction of $\wZ^{0},\wZ^{-1}$ with sector antirenormalization. If the rotation number of $f_n$ is close to $0$, then the $\alpha$ and $\beta$ fixed points are close to each other and $\partial \wZ^0_{f_n}$ travels through the gap between $\alpha= 0$ and $\beta$. The pseudo-Siegel disk $\wZ^{-1}_{f_n}$ is obtained by adding a parabolic fjord to $\wZ^0_{f_n}$. The fjord is bounded by a dam located inside an invariant rectangle, depicted in a pencil style. We say that $\wZ^{0}_{f_n}$ is \emph{regularized} to $\wZ^{-1}_{f_n}$. Under the antirenormalization,  $\wZ^{-1}_{f_n}$ become $\wZ^0_{f_{n-1}}$, and the procedure can be repeated. If the rotation number of $f_n$ is away from $0$, then there is no regularization: $\wZ^{0}_{f_n}=\wZ^{-1}_{f_n}$. 

Informally, since antirenormalization has contracting properties, we expect that the ``non-invariance error'' of $\wZ^{[m]}$ is uniformly bounded by a geometric series. The~\cite{DL_Uniform} justifies such a property by employing the near-degenerate analyses in the original plane of $f_\theta$. The sectorial bounds (cf. Figure~\ref{Fig:SectorRenorm}) are then justified in~\cite{DL:sector-bounds} in the \emph{a posteriori} setting.

The right loop in Figure~\ref{Fig:NearDegRegimIllustr} illustrates a strategy, developed jointly with Lyubich, for extending the theory of pseudo-Siegel disks to the near-neutral setting. On this loop, the antirenormalization is taken rel $\beta$. For many combinatorics, this has been already justified and integrated with \S\ref{ss:NDR:PositiveEntr}-\S\ref{ss:MLC:QL:sat}. If the strategy is fully realized, then this would imply, in particular, the satellite case of the MLC.

\begin{prob}
\label{prob:Satel:unbounded}
Establish pseudo-Siegel a priori bounds in the remaining unbounded satellite ql cases.
\end{prob}

\subsection{Virtual near-Molecule Renormalization.}\label{ss:virtual_renorm} Consider now a general case, when $f_1$ is the first ql renormalization of $f$. The satellite case has been discussed in~\S\ref{ss:MLC:QL:sat} and~\S\ref{ss:MLC:QL:neut}. Let us assume that $f_1$ is primitive ql. The renormalization $f\mapsto f_1$ is encoded by a maximal primitive copy $\MM_1\subset \MM$. There is a canonical sequence of satellite copies of $\MM$

\[\MM\equiv \MM^{(0)} \supsetneq \MM^{(1)}\supsetneq \dots \supsetneq \MM^{(\nn)} \supsetneq \MM^{(\nn+1)},\hspace{2cm} \MM_1\subset \MM\setminus \MM^{(1)}, \qquad \nn\ge 0
\]
such that $\MM^{(m)}$ is a maximal satellite subcopy of $\MM^{(m-1)}$, $\MM_1$ is in the $1/2$-decoration of $\MM^{(\nn+1)}$, and $\MM_1$ is in a small (i.e., not the $1/2$) decoration of $\MM^{(m)}$ for all $m\le \nn$. The copy $\MM_1$ is close to the main Molecule if and only if the period of $\MM^{(\nn+1)}$ is big.

The puzzle levels of~\S\ref{ss:MLC:puzzle} provide a satisfactory control on scales between $\MM^{(\nn)}$ and $\MM_1$ (the latter copy is smaller) and fail on scales between $\MM$ and $\MM^{(\nn)}$, which we refer as \emph{virtual Molecule.} The virtual Molecule case naturally consists of two subcases: \emph{virtual bounded-type satellite ql} (when $\nn\gg 1$ and the relative periods of $\MM^{(s+1)}$ rel $\MM^{(s)}$ are bounded) and \emph{virtual near-Neutral} ($\nn=0$).

\begin{Interprob}\label{prob:Interp}
\label{prob:VirtMolec}
Develop a Virtual Molecule version of the Near-Degenerate Regime.
\end{Interprob}
Problems \ref{prob:Satel:unbounded},~\ref{prob:VirtMolec} and the arguments in~\S\ref{ss:NDR:PositiveEntr},~\S\ref{ss:MLC:puzzle},~\S\ref{ss:MLC:QL:sat}~\S\ref{ss:MLC:QL:neut} should imply the full MLC.

Jointly with Kahn and Lyubich, we put forward a strategy to approach Problem~\ref{prob:VirtMolec} by considering \emph{partially invariant virtual Julia sets} of $\MM^{(s)}$ -- they are connected hulls of the corresponding Cantor small Julia sets within $\Jul_f$ and contain the critical orbit only up to an appropriate number of first returns. \emph{A posteriori}, bounds for virtual Julia sets can be deduced by assuming a uniform hyperbolicity of the renormalization associated with $\MM$; cf. Problem~\ref{conj:hyperb:main}. The strategy towards Problem~\ref{prob:VirtMolec} is to develop such control \emph{a priori}.

\section*{Acknowledgments.} I would like to thank Misha Lyubich, Yusheng Luo, and Kostya Drach for many useful discussions regarding the context of this Note.

\bibliographystyle{amsplain}


\end{document}